\documentclass[12pt]{amsart}
\usepackage{amssymb,amsthm,a4,url}
\usepackage{vmargin}
\newtheorem{thm}{Theorem}
\newtheorem{cor}[thm]{Corollary}
\newtheorem{lem}[thm]{Lemma}
\newtheorem{prop}[thm]{Proposition}
\theoremstyle{definition}

\theoremstyle{remark}
\newtheorem{rem}[thm]{Remark}

\begin{document}
\title[Large balls of torsion and small entropy]
{Infinite groups with large balls of torsion elements and small entropy}
\author{Laurent Bartholdi}
\author{Yves de Cornulier}
\date{October 7, 2005}
\keywords{Grigorchuk group; balls of torsion; algebraic entropy}
\subjclass[2000]{Primary 20F50; Secondary 20F16, 20F05}%
\begin{abstract}
We exhibit infinite, solvable, virtually abelian groups with a
fixed number of generators, having arbitrarily large balls
consisting of torsion elements. We also provide a sequence of
3-generator non-virtually nilpotent polycyclic groups of algebraic
entropy tending to zero. All these examples are obtained by taking
appropriate quotients of finitely presented groups mapping onto
the first Grigorchuk group.
\end{abstract}
\maketitle

The Burnside Problem asks whether a finitely generated group all
of whose elements have finite order must be finite. We are
interested in the following related question: fix $n$ sufficiently
large; given a group $\Gamma$, with a finite symmetric generating
subset $S$ such that every element in the $n$-ball is torsion, is
$\Gamma$ finite? Since the Burnside problem has a negative answer,
{\it a fortiori} the answer to our question is negative in
general. However, it is natural to ask for it in some classes of
finitely generated groups for which the Burnside Problem has a
positive answer, such as linear groups or solvable groups. This
motivates the following proposition, which in particular answers
a question of Breuillard to the authors.

\begin{prop}
For every $n$, there exists a group $G$, generated by a 3-element
subset $S$ consisting of elements of order 2, in which the
$n$-ball consists of torsion elements, and which satisfies one of the
additional assumptions:

(i) $G$ is solvable, virtually abelian, and infinite (more
precisely, it has a free abelian normal subgroup of finite 2-power
index); in particular it is linear.

(ii) $G$ is polycyclic but not virtually nilpotent.

(iii) $G$ is solvable but not polycyclic.\label{prop:main}
\end{prop}

\begin{rem}~\begin{enumerate}
\item The groups in Proposition \ref{prop:main} can actually be
chosen to be 2-generated: indeed, if $G$ is generated by three
involutions $a,b,c$, then the subgroup generated by $ab$ and $bc$
has index at most two.

\item Natural stronger hypotheses are the following: being
linear in fixed dimension; being solvable of given solvability
length. We have no answer in these cases. It is also natural to
ask what happens it we fix a torsion exponent.

\item By~\cite[Corollaire 2, p.~90]{S}, if $G$ is a group and
$S$ is any finite generating subset for which the 2-ball of $G$
consists of torsion elements, then $G$ has Property (FA): every
action of $G$ on a tree has a fixed point. In particular, if $G$
is infinite, then by Stallings' Theorem~\cite{St} it cannot be
virtually free.

\item For every sufficiently large prime $p$, and for all $n$,
there exists a non-elementary, 2-generated word hyperbolic group
in which the $n$-ball consists of elements of
$p$-torsion~\cite{Ol}.

\item We give more precise statements in the sequel: in (i), the
free abelian subgroup can be chosen of index $2^{a_n}$, where
$a_n\sim 13n^k$ (that is, $a_n/(13n^k)\to 1$), where $k\cong 6.60$
is a constant (see Corollary~\ref{cor:n_ball_torsion_sg_normal}).
\end{enumerate}
\label{rem:main_rem}
\end{rem}

With a similar construction, we obtain results on the growth
exponent. Let $G$ be generated by a finite symmetric set $S$, and
denote by $B_n$ the $n$-ball in $G$. Then, by a standard
argument~\cite[Proposition~VI.56]{H}, the limit $h(G,S)=\lim\frac
1n\log(\#(B_n))$ exists. The (algebraic) entropy of $G$ is defined
as $h(G)=\inf_Sh(G,S)$, where $S$ ranges over all finite symmetric
generating subsets of $G$. Osin has proved~\cite{O1,O2} that, for
an elementary amenable finitely generated group, $h(G)=0$ if and
only if $G$ is virtually nilpotent; on the other hand,
Wilson~\cite{W} has constructed a finitely generated group with
$h(G)=0$ which is not virtually nilpotent; see~\cite{B2} for a
simpler example. Relying on former work of Grigorchuk~\cite{G},
Osin observes in~\cite{O2} that there exist elementary amenable
groups (actually they are virtually solvable) with $h>0$ arbitrary
close to 0. This last result can be improved as follows.

\begin{prop}
For every $\varepsilon>0$, there exists a polycyclic, virtually
metabelian, 3-generated group $G$
with~$0<h(G)<\varepsilon$.\label{prop:entropy}
\end{prop}

Propositions~\ref{prop:main} and~\ref{prop:entropy} are obtained by
approximating the Grigorchuk group, first introduced in~\cite{G80}, by
finitely presented groups.

We recall below the definition of a family of 3-generated groups
$\Gamma_n$, which are successive quotients ($\Gamma_{n+1}$ is a
quotient of $\Gamma_n$ for all $n$). These are finitely presented
groups obtained by truncating a presentation of Grigorchuk's
group. They are proved in~\cite{GH} to be virtually direct
products of nonabelian free groups; they have larger and larger
balls of torsion elements, and their entropy tends to zero. We get
Propositions~\ref{prop:main} and~\ref{prop:entropy} by considering
appropriate solvable quotients of the groups~$\Gamma_n$.

Following Lysionok~\cite{L}, the first Grigorchuk group is
presented as follows. We start from the 3-generated group
$\Gamma_{-1}=\langle a,b,c,d\,|\; a^2=b^2=c^2=d^2=bcd=1\rangle$.
Elements $u_n$ and $v_n$ of $\Gamma_{-1}$ are defined below. Then,
for $0\le n\le\infty$, $\Gamma_n$ is defined as the quotient of
$\Gamma_{-1}$ by the relations $u_i$ for $i<n+1$ and $v_i$ for
$i<n$. The first Grigorchuk group $\Gamma=\Gamma_\infty$ has a
wealth of remarkable properties. The most celebrated one is that
$\Gamma$ has non-polynomial subexponential growth~\cite{G}. It is
also a 2-group, i.e.~a group in which every element is of finite
2-power order, and is just-infinite, i.e.\ it is infinite but all its
proper quotients are finite.

We now construct the relations $u_n$ and $v_n$. Consider the
substitution $\sigma$ defined by $\sigma(a)=aca$, $\sigma(b)=d$,
$\sigma(c)=b$, $\sigma(d)=c$; extend its definition to words in
the natural way, and finally observe that it defines a group
endomorphism of $\Gamma_{-1}$. Set $u_0=(ad)^4$, $v_0=(adacac)^4$,
$u_n=\sigma^n(u_0)$, $v_n=\sigma^n(v_0)$.

For all $n\ge -1$, the natural homomorphism
$\Gamma_n\to\mathbf{Z}/2\mathbf{Z}$ sending $b,c,d$ to 0 and $a$
to 1 has kernel $\Xi_n$ of index two.

We will focus on the finitely presented groups $\Gamma_n$ rather than
on $\Gamma$. Individually, and up to commensurability, the structure
of these groups is not of special interest: $\Gamma_n$ is
commensurable to a direct product of $2^n$ non-abelian free
groups~\cite[Proposition 12]{GH}. However, since $\Gamma_\infty$ is
torsion, for all $n$, there exists $i(n)$ such that every element in
the $n$-ball of $\Gamma_{i(n)}$ is torsion.  A quantitative statement
is given in the following proposition, whose proof appears in the
appendix. Let $\lambda\cong 1.25$ be the real root of the polynomial
$2X^3-X^2-X-1$, and set $i(n)=\lfloor\log_\lambda(n)-1\rfloor$.
\begin{prop}
  In the $n$-ball of $\Gamma_{i(n)}$ (for the word metric), every
  element is of $2^{i(n)+1}$-torsion.\label{prop:n_ball_torsion}
\end{prop}

\noindent The following proposition, which
specifies~\cite[Proposition~12]{GH}, describes the structure of $\Gamma_n$.
\begin{prop}
  For every $n\ge 0$, $\Gamma_n$ has a normal subgroup $H_n$ of index
  $2^{\alpha_n}$, where $\alpha_n\le (13\cdot4^n-1)/3$, and $H_n$ is a
  subgroup of index $2^{\beta_n}$ in a finite direct product of $2^n$
  nonabelian free groups of rank 3, where $\beta_n\le (13\cdot
  4^n-15\cdot2^n+2)/3$.\label{prop:normal_product_free}
\end{prop}

\begin{rem}The main difference with~\cite[Proposition 12]{GH} is that
the finite index subgroup they construct is not normal. Of
course one could take a smaller normal subgroup of finite index,
but this one need not \textit{a priori} be of 2-power index, a
fact we require to obtain solvable (and not only virtually
solvable) groups in Propositions~\ref{prop:main}
and~\ref{prop:entropy}.\end{rem}

\noindent We use the following elementary lemma.
\begin{lem}
  Let $G$ be a group, and let $H$ be a proper subgroup of index $2^a$
  in $G$, normalized by a subgroup of index two in $G$. Let $N$ be the
  intersection of all conjugates of $H$. Then $N$ has index $2^b$ in
  $G$, for some integer $b\le 2a-1$.\label{lem:grp_lemma}
\end{lem}
\begin{proof} If $H$ is normal in $G$, the result is trivial.
  Otherwise, consider the unique conjugate $H'\neq H$ of $H$, so that
  $N=H\cap H'$. Taking the quotient by $N$, we can suppose that $H\cap
  H'=\{1\}$ and we are reduced to proving that $G$ is a 2-group of
  order $d\le 2^{2a-1}$. Let $W$ be the normalizer of $H$. Since it
  has index $2$ in $G$, it is normal in $G$, so that $H'\subset W$.
  Since $H$ and $H'$ are both normal subgroups of $W$ and $H\cap
  H'=\{1\}$, $[H,H']=\{1\}$. Accordingly, $HH'$ is a normal subgroup
  of $G$, contained in $W$, and is naturally the direct product of $H$
  and $H'$. The order of $H$ is $d/2^a$, so that the order of $HH'$ is
  $d^2/2^{2a}$, and hence the index of $HH'$ in $G$ is $2^{2a}/d$.
  This proves that $d$ is a power of 2, and $d\le 2^{2a}$; actually
  $d\le 2^{2a-1}$ because $HH'$ is contained in $W$, hence has
  index~$\ge 2$ in $G$.
\end{proof}

\begin{rem}
In Lemma~\ref{lem:grp_lemma}, the assumption that the normalizer
has index at most two is sharp: in the alternating group $A_4$,
there are four subgroups of index 4, all conjugate; they have
pairwise trivial intersection, hence of index 12, which is not a
power of~2.
\end{rem}

Recall that $\Xi_0\subset\Gamma_0$ is a subgroup of index 2; it is
generated by $b$, $c$, $d$, $aba$, $aca$, $ada$.
By~\cite[Proposition 1]{GH}, the assignment $i_0(b)=(a,c)$,
$i_0(c)=(a,d)$, $i_0(d)=(1,b)$ extends to a unique group
homomorphism $i_0:\Xi_0\to\Gamma_{-1}\times\Gamma_{-1}$ such that,
for all $x\in\Gamma_0$, if $i_0(x)=(x_0,x_1)$, then
$i_0(axa)=(x_1,x_0)$. By~\cite[Proposition 10]{GH}, this induces,
for all $n\ge 1$, an injective group homomorphism:
$i_n:\Xi_n\to\Gamma_{n-1}\times\Gamma_{n-1}$.

\begin{proof}[Proof of Proposition~\ref{prop:normal_product_free}] Let
  us proceed by induction on $n$. We start with the essential case
  when $n=0$, worked out in~\cite[Lemma~11]{GH}. Write
  $\Gamma_0=\langle a,b,d|\;a^2=b^2=d^2=(bd)^2=(ad)^4=1\rangle$ (this
  is a Coxeter group). Let $H_0$ be the normal subgroup generated by
  $(ab)^2$. Then, by an immediate verification, $\Gamma_0/H_0$ is
  isomorphic to the direct product of a cyclic group of order 2 and a
  dihedral group of order~8.

We claim that $H_0$ is free of rank 3. Let $L$ be the normal
subgroup of $\Gamma_0$ generated by $ab$: then $L$ has index 4 in
$\Gamma_0$, contains $H_0$, and is shown, in the proof
of~\cite[Lemma~11]{GH}, to be isomorphic to
$\mathbf{Z}\ast(\mathbf{Z}/2\mathbf{Z})$.

Accordingly, by Kurosh's Theorem, if $H_0$ were not free, then it
would contain a conjugate of $ab$, but this is not the case.
Actually $H_0$ is contained in a subgroup of index 8 in
$\Gamma_0$, free of rank 2 (see the proof of~\cite[Lemma~11]{GH}),
hence has rank~3.

Now, for $n\ge 1$, let us suppose that $\Gamma_{n-1}$ has a normal
subgroup $H_{n-1}$ of index $2^{\alpha_{n-1}}$, which embeds as a
subgroup of index $2^{\beta_{n-1}}$ in a direct product of
$2^{n-1}$ non-abelian free groups of rank~3.

The homomorphism $i_n$ described above embeds $\Xi_n$ as a
subgroup of index $8$ in $\Gamma_{n-1}\times\Gamma_{n-1}$. Define
$H'_n=i_n^{-1}(H_{n-1}\times H_{n-1})$: this a normal subgroup of
index $2^k$ in $\Xi_n$, with $k\le 2\alpha_{n-1}$; so $H'_n$ has index
$2^{k+1}$ in $\Gamma_n$.

Then, using Lemma~\ref{lem:grp_lemma}, $H_n=H'_n\cap \,^a\!H'_n$ has index
$2^{\alpha_n}$ in $G$, for some $\alpha_n\le 4\alpha_{n-1}+1$. Combining
the inclusions $H_n\subset H'_n\stackrel{\sim}\to H_{n-1}\times
H_{n-1}\subset {F_3}^{2^{n-1}}\times {F_3}^{2^{n-1}}$, we obtain that
$H_n$ embeds as a subgroup of index $2^{\beta_n}$ in ${F_3}^{2^n}$, with
$\beta_n=[H'_n:H_n]+2[{F_3}^{2^{n-1}}:H_{n-1}]\le 2\alpha_{n-1}+2\beta_{n-1}$.

Define $\alpha'_n=(13\cdot 4^n-1)/3$ and
$\beta'_n=(13\cdot 4^n-15\cdot 2^n+2)/3$. Then $\alpha'_0=4$,
$\beta'_0=0$, and they satisfy, for all $n$:
$\alpha'_n=4\alpha'_{n-1}+1$ and
$\beta'_n=2\alpha'_{n-1}+2\beta'_{n-1}$. Therefore, an immediate
induction gives $\alpha_n\le\alpha'_n$ and $\beta_n\le\beta'_n$
for all~$n\ge0$.\end{proof}

It is maybe worthwhile to restate the result avoiding reference to the
particular sequence~$\Gamma_n$.

\begin{cor}\label{cor:epi_to_torsion}
  For every finitely presented group $G$ mapping onto the first
  Grigorchuk group $\Gamma$, there exist normal subgroups
  $N\subsetneq H$ in $G$, with $H$ of finite 2-power index,
  such that $H/N$ is isomorphic to a finite index subgroup in a direct
  product of non-abelian free groups.
\end{cor}
\begin{proof} Let $p:G\to\Gamma$ be onto. Since $G$ is finitely
  presented, $p$ factors through $\Gamma_n$ for some sufficiently
  large $n$, so that there exists a map $p'$ of $G$ onto $\Gamma_n$.
  Then take $N=\text{Ker}(p')$ and $H=p'^{-1}(H_n)$.
\end{proof}

\noindent Combining Propositions~\ref{prop:n_ball_torsion}
and~\ref{prop:normal_product_free}, we also obtain the following
statement:
\begin{cor}
  In the group $\Gamma_{i(n)}$, the $n$-ball consists of
  $2^{i(n)+1}$-torsion elements, and there exists a normal subgroup of
  index $2^{\alpha_{i(n)}}$, which embeds in a direct product of free
  groups, with $\alpha_{i(n)}\le (13\cdot n^{\log_\lambda(4)}-1)/3$,
  and $\log_\lambda(4)\cong 6.60$.\qed\label{cor:n_ball_torsion_sg_normal}
\end{cor}

\begin{proof}[Proof of Proposition~\ref{prop:main}] By
  Proposition~\ref{prop:n_ball_torsion}, we may take $i$ sufficiently
  large so that the $n$-ball of $\Gamma_i$ consists of torsion
  elements. Since $H_i$ is a finite index subgroup in a nontrivial
  direct product of free groups (see
  Proposition~\ref{prop:normal_product_free}), it has infinite
  abelianization.  There is a short exact sequence $$1\to
  H_i/[H_i,H_i]\to \Gamma_i/[H_i,H_i]\to \Gamma_i/H_i\to 1\,.$$
  Accordingly, $G=\Gamma_i/[H_i,H_i]$ is an infinite, virtually
  abelian group, in which the $n$-ball consists of torsion elements.
  Moreover, since $\Gamma_i/H_i$ is a finite 2-group, $G$ is also
  solvable. This proves (i).

For (iii), take, instead, $G=\Gamma_i/[[H_i,H_i],[H_i,H_i]]$.
Since $H_i$ maps onto a non-abelian free group, its
metabelianization is not virtually polycyclic, so that $G$ is
virtually metabelian, but not virtually polycyclic.

For (ii), take a morphism of $H_i$ onto a polycyclic group $W$
which is not virtually nilpotent, and let $K$ be the kernel of
this morphism. Since the normalizer of $K$ has finite index in
$\Gamma_i$, $K$ has finitely many conjugates $K_1,\dots,K_\ell$.
Set $L=\bigcap_{j=1}^\ell K_i$. Then the diagonal map
$H_i/L\to\prod_{j=1}^\ell H_i/K_i$ is injective, hence embeds
$H_i/L$ in $W^\ell$. On the other hand, observe that $H_i/L$
projects onto $W$, so is not virtually nilpotent. It follows that
$G=\Gamma_i/L$ is polycyclic but not virtually nilpotent. If $W$
has been chosen metabelian, then we also have that $G$ is
virtually metabelian.
\end{proof}

\begin{proof}[Proof of Proposition~\ref{prop:entropy}] Keep
the last construction in the previous proof. Then
$h(\Gamma_i/L)\le h(\Gamma_i)$. Moreover, $h(\Gamma_i/L)>0$ since
$\Gamma_i/L$ is solvable but not virtually nilpotent~\cite{O1}. On
the other hand, it is proved in~\cite{GH} that $h(\Gamma_i)\to 0$.
Thus we can obtain $h(G)$ arbitrarily small.\end{proof}

\begin{rem} Consider for every $i$ an infinite quotient $Q_i$ of
  $\Gamma_i$. In the topology of marked groups (defined in~\cite{G};
  see also, for instance,~\cite{Cha}), the sequence $(Q_i)$ converges
  to the Grigorchuk group $\Gamma$. Indeed, otherwise by compactness
  it would have another cluster point, which would be a proper
  quotient of $\Gamma$, and therefore would be finite. This is a
  contradiction since the infinite groups form a closed subset in the
  topology of marked groups.
\end{rem}

\appendix

\section*{Appendix}
We gather here the technical results concerning torsion in the groups
$\Gamma_n$. They are slight modifications of results in the
papers~\cite{GH,B}.

Recall that $\lambda\cong 1.25$ denotes the real root of the
polynomial $2X^3-X^2-X-1$.  We introduce on $\Gamma_{-1}$ (hence on
all its quotients) the metric $|\cdot|$ of~\cite{B}: it is defined by
attributing a suitable weight to each of the generators $a$, $b$, $c$,
$d$: $|a|=2(\lambda-1)\cong 0.47$, $|b|=1-|a|=\lambda^{-3}\cong 0.53$,
$|c|=2\lambda^2-3\lambda+1\cong 0.34$, and
$|d|=-2\lambda^2+\lambda+2\cong 0.19$. {\bf Throughout this appendix,
  unless explicitly stated, the balls and the lengths are meant in the
  sense of this weighted metric.}

To check that the length of $a$, $b$, $c$, $d$ is exactly the
weight we have imposed, it suffices to check this in the
abelianization of $\Gamma_{-1}$, the $\mathbf{F}_2$-vector space
with basis $(a,b,d)$ (which is also the abelianization of all
$\Gamma_n$). There, it is a straightforward verification that the
mapping $|\cdot|$ just defined extends to a length function by
setting $|a\xi|=|a|+|\xi|$ for all $\xi\in\{b,c,d\}$.

Observe that if $\xi\in\{b,c,d\}$, and $i_n(\xi)=(\xi_0,\xi_1)$,
we have
\begin{equation} |\xi_0|+|\xi_1|=\lambda^{-1}(|\xi|+|a|).
\label{e:eeee}\end{equation}

\begin{lem}
Let $x\in\Gamma_0$ be any element. Set $x'=x$ if $x\in\Xi_0$ and
$x'=xa$ otherwise; and set $i_0(x')=(x_0,x_1)$.

Then $|x_0|+|x_1|\le\lambda^{-1}(|x|+|a|)$.

Suppose moreover that $x$ is of minimal length among its
conjugates, and that $x\notin\{b,c,d\}$. Then
$|x_0|+|x_1|\le\lambda^{-1}|x|$.\label{lem:ineq_x0+x1}
\end{lem}

\begin{proof} Fix $x\in\Gamma_0$, and let $w$ be a word in the letters
  $\{a,b,c,d\}$, of minimal length\footnote{If $w=u_1\dots u_n$, the
    length of $w$ is defined as $|u_1|+\dots+|u_n|$.}, representing
  $x$. Since every element in $\{b,c,d\}$ is the product of the two
  others, $w$ can be chosen so that no two consecutive letters are
  in~$\{b,c,d\}$.

Suppose now that $x$ is of minimal length within its conjugacy
class and that $w$ is not a single letter. Maybe conjugating $x$
by the last letter of $w$, we can suppose that $w$ ends with the
letter~$a$. The minimality assumption then implies that $w$ begins
with a letter in $\{b,c,d\}$.

\begin{description}
    \item[First case] $x\in\Xi_0$. Write
    $w=\xi^1(a\xi^2a)\dots\xi^{2n-1}(a\xi^{2n}a)$, where
    $\xi^i\in\{b,c,d\}$ for $i=1,\dots,2n$. Write $i(x)=(x_0,x_1)$ and
    $i(\xi^i)=(\xi^i_0,\xi^i_1)$, so that
    $i(\xi^i)=(\xi^i_1,\xi^i_0)$. Then

\begin{align*}
  |x_0|+|x_1|&\le (|\xi^1_0|+|\xi^2_1|+\dots+|\xi^{2n-1}_0|+|\xi^{2n}_1|)+
(|\xi^1_1|+|\xi^2_0|+\dots+|\xi^{2n-1}_1|+|\xi^{2n}_0|)\\
&=(|\xi^1_0|+|\xi^1_1|)+(|\xi^2_0|+|\xi^2_1|)+\dots+(|\xi^{2n}_0|+|\xi^{2n}_1|).
\end{align*}

    By (\ref{e:eeee}), we get
$$|x_0|+|x_1|\le\lambda^{-1}\sum_{i=1}^{2n}(|\xi^i|+|a|)$$
    On the other hand, $|x|=\sum_{i=1}^{2n}(|\xi^i|+|a|)$, so that
    finally $|x_0|+|x_1|\le\lambda^{-1}|x|$.

\item[Second case] $x\notin\Xi_0$, so that $xa\in\Xi_0$.
    Write $w=\xi^1(a\xi^2a)\dots\xi^{2n-1}(a\xi^{2n}a)\xi^{2n+1}a$, so
    that $\xi^1(a\xi^2a)\dots\xi^{2n-1}(a\xi^{2n}a)\xi^{2n+1}$
    represents $xa$ in $\Gamma_{0}$. Write $i_0(xa)=(x_0,x_1)$, and
    $i_0(\xi^i)=(\xi^i_0,\xi^i_1)$. Then

\begin{align*}
  |x_0|+|x_1|&\le
  (|\xi^1_0|+|\xi^2_1|+\dots+|\xi^{2n-1}_0|+|\xi^{2n}_1|+|\xi^{2n+1}_0|)\\
  &\kern+2cm{}+ (|\xi^1_1|+|\xi^2_0|+\dots+|\xi^{2n-1}_1|+|\xi^{2n}_0|+|\xi^{2n+1}_1|)\\
  &=(|\xi^1_0|+|\xi^1_1|)+(|\xi^2_0|+|\xi^2_1|)+\dots+(|\xi^{2n+1}_0|+|\xi^{2n+1}_1|)\\
  &=\lambda^{-1}\bigg(\sum_{i=1}^{2n+1}(|\xi^i|+|a|)\bigg)\quad\text{again by (\ref{e:eeee})}.
\end{align*}

    Since $|x|=\sum_{i=1}^{2n+1}(|\xi^i|+|a|)$, we get
    $|x_0|+|x_1|\le\lambda^{-1}|x|$.
  \end{description}

The other inequality $|x_0|+|x_1|\le\lambda^{-1}(|x|+|a|)$ is
proved similarly: we must deal with the following
cases:\begin{itemize}
 \item $w$ begins and ends with the letter $a$: considering whether
 or not
 $x\in\Xi_0$, in both cases we obtain the stronger inequality
 $|x_0|+|x_1|\le\lambda^{-1}(|x|-|a|)$.
 \item $w$ begins and ends with letters in $\{b,c,d\}$: considering
 whether or not
 $x\in\Xi_0$, in both cases we obtain the inequality
 $|x_0|+|x_1|\le\lambda^{-1}(|x|+|a|)$.
 \item $w$ begins with the letter $a$ and ends with a letter in
$\{b,c,d\}$: in this case, replacing $x$ by $x^{-1}$ and $w$ by
$w^{-1}$ (this is the word $w$ read from right to left
---~recall that the generators are involutions), we reduce to the
case, already carried out, in which $w$ begins with a letter in
$\{b,c,d\}$ and ends with the letter~$a$, obtaining the inequality
$|x_0|+|x_1|\le\lambda^{-1}|x|$.
\end{itemize}
Since the verifications are similar to the computations above and
since we do not use this case in the sequel, we omit the
details.\end{proof}

\begin{lem}
For every $n\ge -1$, and every element in the open
$\lambda^{n-1}$-ball of $\Gamma_{-1}$, its image in $\Gamma_n$ is
of $2^{n+1}$-torsion.\label{lem:torsion_ball_Gamma_n}
\end{lem}

\begin{proof} For $n=-1$, $\lambda^{-2}=|a|+|d|\cong 0.66$, and the elements
in the open $\lambda^{-2}$-ball are $1$, $a$, $b$, $c$, and $d$,
and are of 2-torsion in~$\Gamma_{-1}$.

For $n=0$, $\lambda^{-1}=|a|+|c|\cong 0.81$, and the elements in
the open $\lambda^{-1}$-ball are, besides the elements in the open
$\lambda^{-2}$-ball already quoted, $ad$ and its inverse $da$,
which are of 4-torsion in~$\Gamma_0$.

We can start an induction, and suppose that, for some $n\ge 1$, we
have already proved that for every element in the open
$\lambda^{n-2}$-ball of $\Gamma_{-1}$, its image in $\Gamma_{n-1}$
is of $2^{n}$-torsion. Pick $x$ in the open $\lambda^{n-1}$-ball
of $\Gamma_{-1}$. We want to show that $x^{2^{n+1}}=1$. We can
suppose that $x$ is of minimal length among its conjugates, and
that $x\notin\{b,c,d\}$. Define $x'$ as in
Lemma~\ref{lem:ineq_x0+x1}, i.e.~$\{x'\}=\{x,xa\}\cap\Xi_{-1}$.
Denote by $[\cdot]$ the projection of $\Gamma_{-1}$ onto
$\Gamma_{0}$. Set $i_0([x'])=(x_0,x_1)$.

\begin{description}
\item[First case] $x\in\Xi_{-1}$, i.e.~$x=x'$. By Lemma
   ~\ref{lem:ineq_x0+x1}, for $i=0,1$, $|x_i|\le
    |x_0|+|x_1|\le\lambda^{-1}|x|\le\lambda^{n-2}$. By the induction
    hypothesis, $x_0$ and $x_1$ are of $2^{n}$-torsion in
    $\Gamma_{n-1}$. Since $i_0$ induces an injection of $\Xi_n$ into
    $\Gamma_{n-1}\times\Gamma_{n-1}$, this implies that $x$ is of
    $2^n$-torsion, hence of $2^{n+1}$-torsion in~$\Gamma_n$.

\item[Second case] $x\notin\Xi_{-1}$, i.e.~$x'=xa$. Then
    $x^2\in\Xi_{-1}$, and
$$i_0([x^2])=i_0([xaaxaa])=i_0([xa])i_0([a(xa)a])=(x_0x_1,x_1x_0),$$
    which is conjugate to $(x_0x_1,x_0x_1)$ in
    $\Gamma_{-1}\times\Gamma_{-1}$. By Lemma~\ref{lem:ineq_x0+x1}, we
    have $|x_0x_1|\le|x_0|+|x_1|\le\lambda^{-1}|x|\le\lambda^{n-2}$. By
    the induction hypothesis, $x_0x_1$ is of $2^n$-torsion in
    $\Gamma_{n-1}$, so $i_0([x^2])$ is of $2^n$-torsion in
    $\Gamma_{n-1}\times\Gamma_{n-1}$. Since $i_0$ induces an injection
    of $\Xi_n$ into $\Gamma_{n-1}\times\Gamma_{n-1}$, this implies
    that $x^2$ is of $2^n$-torsion in $\Xi_n\subset\Gamma_n$, hence
    $x$ is of $2^{n+1}$-torsion in $\Gamma_n$.
  \end{description}
\end{proof}

\begin{proof}[Proof of Proposition~\ref{prop:n_ball_torsion}] Suppose
  that $x$ has word length $\le n$. Then
  \[|x|\le
  |b|n=\lambda^{-3}n=\lambda^{\log_\lambda(n)-3}<\lambda^{i(n)-1}.\]
  By Lemma~\ref{lem:torsion_ball_Gamma_n}, the image of $x$ in
  $\Gamma_{i(n)}$ is of $2^{i(n)+1}$-torsion.
\end{proof}

\bigskip

\noindent\textbf{Acknowledgement}. We thank Emmanuel Breuillard for
encouragement and valuable discussions. We also thank Goulnara
Arzhantseva and Pierre de la Harpe for useful remarks.

\bigskip
\footnotesize
\noindent Laurent Bartholdi\\
E-mail: \url{laurent.bartholdi@epfl.ch}\\
\'Ecole Polytechnique F\'ed\'erale de Lausanne (EPFL)\\
Institut de Math\'ematiques B (IMB)\\
CH-1015 Lausanne, Switzerland

\medskip

\noindent Yves de Cornulier\\
E-mail: \url{decornul@clipper.ens.fr}\\
Universit\'e de Neuch\^atel\\
Institut de Math\'ematiques\\
Rue \'Emile Argand 11\\
CH-2007 Neuch\^atel, Switzerland


\begin{thebibliography}{KM98b}

\bibitem[B1]{B} Laurent {\sc Bartholdi}. \newblock {\em The growth of
    Grigorchuk's torsion group}. \newblock Int. Math. Res. Not. {\bf
    20}, 1049-1054 (1998).

\bibitem[B2]{B2} Laurent {\sc Bartholdi}.
\newblock {\em A Wilson group of non-uniformly exponential
growth}. C. R. Math. Acad. Sci. Paris {\bf 336}, no.7, 549-554
(2003).

\bibitem[C]{Cha} Christophe {\sc Champetier}.
\newblock {\em L'espace des groupes de type fini}, Topology~\textbf{39}, no.4 657-680 (2000).

\bibitem[G1]{G80} Rostislav {\sc Grigorchuk}. \newblock {\em On the Burnside problem for periodic
groups}. \newblock Funct. Anal. Appl. {\bf 14}, 41-43, 1980
(Russian original: Funktsional. Anal. i Prilozhen {\bf 14}, 53-54,
1980).

\bibitem[G2]{G} Rostislav I. {\sc Grigorchuk}. \newblock {\em Degrees of growth
of finitely generated groups, and the theory of invariant means}.
\newblock Math. USSR Izv. {\bf 25}:2, 259-300, 1985 (Russian original: Izv. Akad Nauk. SSSR
Ser. Mat. {\bf 48}, no.5, 939-985, 1984).

\bibitem[GH]{GH} Rostislav I. {\sc Grigorchuk}, Pierre {\sc de la Harpe}. \newblock
{\em Limit behaviour of exponential growth rates for finitely
generated groups}. \newblock Monogr. Enseign. Math. {\bf 38},
351-370 (2001).

\bibitem[H]{H}
Pierre~{\sc de la Harpe}. \newblock ``Topics in geometric group
  theory''. \newblock University of Chicago Press, Chicago, IL, 2000.

\bibitem[L]{L} Igor G. {\sc Lysionok}. \newblock {\em A set of defining relations for the Grigorchuk
group}. \newblock Math. Notes Acad. Sc. USSR {\bf 38}:4, 784-792,
1985 (Russian original: Mat. Zametki {\bf 38}, 503-511, 1985).

\bibitem[Ol]{Ol} Alexander {\sc Olshanskii}. \newblock {\em An
infinite group with subgroups of prime orders}. \newblock Math.
USSR Izv. {\bf 16}, 279-289, 1981 (Russian original: Izv. Akad.
Nauk SSSR Ser. Mat. {\bf 44}, 309-321, 1980).

\bibitem[Os1]{O1} Denis V. {\sc Osin}. \newblock {\em The entropy of solvable groups}.
\newblock Ergod. Theory \& Dynam. Sys. {\bf 23}, no.3, 907-918 (2003).

\bibitem[Os2]{O2} Denis V. {\sc Osin}. \newblock {\em Algebraic entropy of elementary amenable groups}.
\newblock Preprint 2005, to appear in Geom. Dedicata.

\bibitem[Se]{S} Jean-Pierre {\sc Serre}. \newblock ``Arbres,
amalgames, $\textnormal{SL}_2$''. \newblock Ast\'erisque {\bf 46},
SMF, 1977.

\bibitem[St]{St} John R. {\sc Stallings}. \newblock {\em On torsion-free
groups with infinitely many ends}. \newblock Ann. of Math. {\bf 88},
312-334 (1968).

\bibitem[W]{W} John S. {\sc Wilson}. \newblock {\em On exponential
growth and uniformly exponential growth for groups}. Invent. Math.
{\bf 155}, 287-303 (2004).
\end{thebibliography}
\end{document}